\documentclass[11pt]{article}

\usepackage[centertags]{amsmath}
 \usepackage{amsfonts}
\usepackage{amssymb}
\usepackage{makeidx}
\usepackage{newlfont}
\makeindex 
\usepackage{amsthm} 
\usepackage{stmaryrd}
\usepackage[]{titletoc}  
\usepackage{mathrsfs}
\usepackage{stmaryrd}
 \usepackage{graphicx}  
 \usepackage{xypic}
\usepackage{hyperref}

\newlength{\defbaselineskip}
\newcommand{\setlinespacing}[1]%
           {\setlength{\baselineskip}{#1 \defbaselineskip}}

\theoremstyle{plain}
\newtheorem{thm}{Theorem}[section]
\newtheorem{cor}[thm]{Corollary}
\newtheorem{lem}[thm]{Lemma}

\newtheorem{exam}[thm]{Example}
\newtheorem{rem}[thm]{Remark}

\newcommand{\hard}{H^2(\mathbb{D})}

\makeatletter\@addtoreset{equation}{section} \makeatother
\begin{document}
\title {Cowen's class and Thomson's class}
\author{Kunyu Guo\ \  \  Hansong Huang }
\date{}
 \maketitle \noindent\textbf{Abstract:}
In studying commutants of analytic Toeplitz operators,   Thomson
 \cite{T2} proved  a remarkable theorem which states that under a mild
condition, the commutant  of an analytic Toeplitz operator is equal to that of   Toeplitz operator defined by a finite Blaschke product.
Cowen\cite{Cow1} gave an significant improvement of Thosom's
result. In this paper, we will  present examples in Cowen's class which does not lie in Thomson's class.


\section{Introduction}
~~~~In this paper, $\mathbb{D}$ denotes the unit disk in the complex plane $\mathbb{C}$, and
 $\hard$ denotes the Hardy space on $\mathbb{D}$, which consists of all holomorphic functions   whose Taylor coefficients at $0$
  are  square summable. Let $H^\infty(\mathbb{D})$
  denote  the set of all bounded holomorphic functions over $\mathbb{D}$. For each function $h $ in  $ H^\infty(\mathbb{D})$,
  $T_h$ denotes the Toeplitz operator  on the Hardy space $\hard$, that is, the multiplication operator defined by the symbol $h$.

In \cite{DW}  Deddens and  Wong  raised several questions about the
commutants for  analytic  Toeplitz operators defined on the Hardy
space  $H^2(\mathbb{D})$.
 One of them asks that for a function $\phi\in  H^\infty(\mathbb{D})$, whether there is
  an inner function $\psi$ such that $\{T_\phi\}'=\{T_\psi\}'$ and that
  $\phi=h\circ \psi$  for some $h\in  H^\infty(\mathbb{D})$, where $\{T_\phi\}'=\{A\in B(\hard):\, AT_\phi=T_\phi A\} $ is the commutant of $T_\phi$.
     Baker, Deddens and  Ullman\cite{BDU} proved that for an entire function $f$,
 there is a positive integer $k$ such that $\{T_f\}'=\{T_{z^k}\}'$. Later, by using function-theoretic techniques,
 Thomson   gave the following remarkable result, see \cite{T2}. 
\begin{thm}[Thomson]
  Suppose that $\phi\in  H^\infty(\mathbb{D})$, and there exist uncountably many points
$\lambda $ in $\mathbb{D}$ such that the inner part  of
$\phi-\phi(\lambda)$ are finite Blaschke products. Then there exists
a finite Blaschke product $B$  and an $H^\infty$-function  $\psi$
such that $\phi =\psi(B)$ and $\{T_\phi\}'= \{T_B\}'$ holds on the Hardy space $H^2(\mathbb{D})$.
\label{Tmth1}
\end{thm}
As a consequence,   the following is immediate, see \cite{T1}.
\begin{cor}
Let $\phi$ be a nonconstant  function that is holomorphic  on the
closed unit disk $\overline{\mathbb{D}}$. Then there exists a finite
Blaschke product $B$ and a \linebreak $\psi\in $
$H^\infty(\overline{\mathbb{D}})$ such that $\phi =\psi(B)$ and
$\{T_\phi\}'= \{T_B\}'$ holds on $\hard$. In particular, if
$\phi$ is entire, then $\psi$ is entire and $B(z)=z^n$ for some
positive integer $n$. \label{Tmth2}
\end{cor}
 One related topic is to study $\{T_h, T_h^*\}'$ for $h\in H^\infty(\mathbb{D})$. Now  put \linebreak $\mathcal{V}^*(h)=\{T_h, T_h^*\}'$,
which turns out to be a von Neumann algebra. It is interesting because for each $H^\infty$-function $\phi$
satisfying the condition in Theorem \ref{Tmth1}, there is
 necessarily a finite Blaschke product $B$ satisfying $\mathcal{V}^*(\phi)=\mathcal{V}^*(B).$
 It is worthwhile to mention that Theorem 1.1  holds not only on the Hardy space,
  but also on the Bergman space. On the Hardy space, the structure of $\mathcal{V}^*(B) $ is clear. However, on the Bergman space it is not easy.
   It is known that $\mathcal{V}^*(B) $ is always nontrivial\cite{HSXY}, and very recently, it is shown that $\mathcal{V}^*(B) $
    is abelian for each finite Blaschke product\cite{DPW}.  Along this line, there are a lot of work,  also refer to
     \cite{Cow1}-\cite{CW},\cite{T1}-\cite{T4},\cite{SW,Zhu1,GSZZ,GH1,GH3,Sun,SZZ1,SZZ2}.

Later, Cowen\cite{Cow1} gave an   significant  extension of Theorem 1.1,
as follows.
\begin{thm}[Cowen]
  Suppose that $\phi\in  H^\infty(\mathbb{D})$, and there exist a point
$\lambda $ in $\mathbb{D}$ such that the inner part  of
$\phi-\phi(\lambda)$ is a finite Blaschke product. Then there exists \label{t13} a finite Blaschke product $B$  and an
$H^\infty$-function $\psi$ such that $\phi =\psi(B)$ and
$\{T_\phi\}'= \{T_B\}'$ holds on the Hardy space $H^2(\mathbb{D})$.
\end{thm}
 However, as Cowen proved Theorem \ref{t13}\cite{Cow1}, he did  not know whether there exists any function $f$ for which

 \vskip2mm $\{ \lambda \in \mathbb{D}|$ \emph{the inner
part of $f - f(\lambda)$ is a finite Blaschke product}$\}$
 \vskip2mm  \noindent is a nonempty
countable set. If there were no functions with this property, then  every
function which satisfies the hypotheses of  Theorem \ref{t13} would also
satisfy the hypotheses of Thomson's theorem.
  This paper is to  construct a function satisfying the property as mentioned above, and thus
    explains why Cowen's generalization is essential.

Before continuing,  we present two conditions on functions in
$H^\infty(\mathbb{D})$. It is convenient to call the assumption in
Theorem  \ref{t13} \emph{Cowen's condition}. That is,  if $h$ is a
function  in $H^\infty(\mathbb{D})$ such that  for some $\lambda $
in $\mathbb{D}$ the inner part of \linebreak $h-h(\lambda)$ is a finite
Blaschke product, then $h$ is said to satisfy \emph{Cowen's
condition}. Similarly, for a function $\phi$ in
$H^\infty(\mathbb{D})$, if there are uncountably many $\lambda $ in
$\mathbb{D}$ such that the inner part of $\phi-\phi(\lambda)$ is a
finite Blaschke product, then $\phi$ is said to satisfy
\emph{Thomson's condition}.  All functions satisfying Thomson's
condition consist of a set, called \emph{Thomson's class},
 and the set of all  functions of  satisfying  Cowen's condition, is  called \emph{Cowen's class}.
Clearly, Thomson's class is contained in Cowen's class.

 One natural question is whether these two classes are the same.
  Cowen \cite{Cow1} raised    it as a question precisely as follows: \vskip2mm

\emph{Is there  a function in Cowen's class which does not lie in
Thomson's class?} \vskip2mm
\noindent  This is  to ask whether Thomson's class is  properly contained in  Cowen's class.
  The  following example provides an affirmative answer, and the
  details will be given in Section 2. 
\begin{exam} Pick $a\in \mathbb{D}-(-1,1)$, and denote by $B$ the
 thin Blaschke product with only simple zeros: $a$ and $1-\frac{1}{n!}(n\geq 2)$.
 By Riemann  mapping theorem, there is a conformal map $h$ from the unit disk onto $\mathbb{D}-[0,1)$.
Take a such $h$, define $\phi=  B\circ h $ and put
$b=h^{-1}(a).$ Later  one
  will see that  the inner part of $\phi-\phi(b)$ is a finite Blaschke product;
  and for any $\lambda\in \mathbb{D}-\{b\}$, the inner part of $\phi-\phi(\lambda)$ is
never a finite Blaschke product. \label{ex14} Thus, $\phi$ is not in Thomson's class,   though it lies in Cowen's class.
\end{exam}

 \section{The construction of examples}
  In this  section
we will provide  examples in  Cowen's class, but not lying in Thomson's class.

Below,   $d$  will denotes the hyperbolic metric on
$\mathbb{D}$; that is,
$$d(z,w)=\big|\frac{z-w}{1- \overline{z}w }\big| ,\ z,w\in \mathbb{D}.$$First, recall that a   Blaschke product with the zero sequence $\{z_k\}$
is called a \emph{thin Blaschke product} if $$\lim_{k\to\infty}\prod_{j:j\neq k} d(z_k,z_j) =1.$$
The following example\cite{GH2} presents the construction for a class of
thin Blaschke products.
\begin{exam} This example comes from \cite{GH2}. For the reader's convenience, we present its details.

 First let us make an observation from \cite[pp. 203,204]{Hof}: 
if $\{w_n\}$ is  a   sequence  in   $\mathbb{D}$ satisfying
\begin{equation} \frac{1-|w_n|}{1-|w_{n-1}|}\leq c<1,\label{5341} \end{equation}
then \begin{equation} \prod_{j:j\neq k} d(w_k,w_j) \geq
\big(\prod_{j=1}^\infty \frac{1-c^j}{1+c^j}\big)^2.  \label{5342}
\end{equation} Notice that the right hand side tends to $1$ as $c$
tends to $ 0$.

Let $\{c_n\}$ be a  sequence satisfying $c_n>0$ and $
\lim\limits_{n\to \infty}c_n=0.$ Suppose that $\{z_n\}$ is a
sequence of points in the open unit disk $\mathbb{D}$ such that
$$ \frac{1-|z_n|}{1-|z_{n-1}|}=c_n.$$
We will show that $\{z_n\}$ is a thin Blaschke sequence. To see
this,  given a positive integer $m$, let $k>m$ and consider the
product
\begin{equation} \prod_{j:j\neq k} d(z_k,z_j) \equiv \prod_{1\leq j\leq m} d(z_k,z_j) \prod_{j>m,j\neq k} d(z_k,z_j). \label{v.25} \end{equation}
Now write  $$d_m =\sup \{c_j:j\geq m\}.$$
 Since $ \lim\limits_{n\to \infty}c_n=0,$ then $ \lim\limits_{m\to \infty}d_m=0.$
For any $k>m$,  $$ \frac{1-|z_k|}{1-|z_{k-1}|}=c_k\leq d_m,$$ and
thus
$$  \prod_{j>m,j\neq k} d(z_k,z_j)\geq \big(\prod_{j=1}^\infty \frac{1-d_m^j}{1+d_m^j}\big)^2, \, k>m, $$
which implies that $$    \prod_{j>m,j\neq k} d(z_k,z_j) $$ tends to
$1$ as $m\to \infty$. For any $\varepsilon >0$, there is an $m_0$
such that
$$    \prod_{j>m_0,j\neq k} d(z_k,z_j)\geq 1- \varepsilon .$$
On the other hand,  $$\lim\limits_{k\to \infty } \prod_{1\leq j\leq
m_0} d(z_k,z_j)=1,$$ which, combined with (\ref{v.25}), implies that
$$ \liminf_{k\to\infty} \prod_{j:j\neq k} d(z_k,z_j) \geq
1-\varepsilon .$$ Therefore, by the arbitrariness of $\varepsilon$,
$$\lim_{k\to\infty}\prod_{j:j\neq k} d(z_k,z_j) =1,$$
which implies that $\{z_k\} $ is a thin Blaschke product. \label{21}
  For example, $\{z_n\}$ is a thin
Blaschk sequence if we put $|z_n|=1- \frac{1}{n!}.$
\end{exam}

  The following is  the restatement  
 of Example \ref{ex14}.
\begin{thm}  There is a holomorphic function $\phi$ from $\mathbb{D}$ onto $\mathbb{D}$, such that
the inner part of $\phi-w\ (w\in \mathbb{D})$ is a finite Blaschke
product  if and only if $w=0.$ \label{351}
\end{thm}
\begin{cor}  There is a holomorphic function $\phi$ from $\mathbb{D}$ onto $\mathbb{D}$, whose zero set $Z(\phi)$ is
a nonempty finite subset of $\mathbb{D}$, such that
the inner part of $\phi-\phi(\lambda)\ (\lambda\in \mathbb{D})$ is a finite Blaschke
product  if and only if    $\lambda  $ lies in $Z(\phi)$.
\end{cor}
\begin{rem}  Corollary 2.3 immediately shows that Thomson's class is properly contained in Cowen's class.
\end{rem}
  This section mainly furnishes the details for Theorem \ref{351}.
    In fact, to prove Theorem \ref{351}, let us verify the details of Example \ref{ex14}.

    Remind that
  $a$ is a point in $ \mathbb{D}-\mathbb{R}$, and $B$  denotes the
 Blaschke product with only simple zeros: $a$ and $1-\frac{1}{n!}(n\geq 2)$.
By Example   \ref{21}, this  Blaschke product $B$ is a thin
 Blaschke product.

   We will construct a conformal map $g$ from  $\mathbb{D}-[0,1)$ onto $\mathbb{D}$
and set $h=g^{-1}$. Then define $\phi= B\circ h$, which turns out to
be the function as in Theorem \ref{351}.
\begin{lem}With $\phi$ defined as above, for any $w\in \mathbb{D}-\{0\}$, the inner part of
$\phi-w$ \label{00l}
is never a finite Blaschke product. \end{lem}  \begin{proof} Remind
that $\phi=B\circ h,$ where $h$ is a conformal map from $\mathbb{D}$
onto $\mathbb{D}-[0,1).$ To prove Lemma  \ref{00l} it suffices to show that
$ B|_{\mathbb{D}-[0,1)}$ attains each nonzero value $w$ in
$\mathbb{D}$ for infinitely many times.
By   \cite[Lemma 3.2(3)]{GM}
 each value in $\mathbb{D}$ can be achieved for infinitely many times by $B$, and then it is enough to show that
 $B|_{[0,1)}$  attains each nonzero value  $w$ in $B\big([0,1)\big)$ for  finitely many times.

For this, notice that  $\varphi_a$ maps $[-1,1]$ to a circular arc
in $\overline{\mathbb{D}} .$ Then
 for each fixed $r\in [-1,1]$,     the argument function $ \mathrm{arg}\  \varphi_a(t)|_{[-1,1]}$ of $\varphi_a(t)$
  attains the value $ \mathrm{arg}\ \varphi_a(r)$ for at most $k_0$ times
(say, $k_0=2$). Here, the value of $\mathrm{arg}$ is required to be
in $[0,2\pi)$. Write $$B=\varphi_a B_0 ,$$ where $B_0$ is a Blaschke
product,
 and it is clear that
 $B_0(r)\in (-1,1) .$
 If $B_0(r)\neq 0,$ then
 $$\mathrm{either}\ \  \mathrm{arg} \ B(r)=\mathrm{arg} \ \varphi_a (r)\
 \ \mathrm{or} \ \ \mathrm{arg} \ B(r)=\mathrm{arg} (\varphi_a (r))+\pi \mod 2\pi .$$
 Therefore, for each $r\in [0,1)-Z(B)$,
$ \mathrm{arg} \ B$ attains the value  $\mathrm{arg}\  B(r)$ for at
most $2k_0$ times on $[0,1)-Z(B)$. Thus  $B|_{[0,1)}$  attains each
nonzero value  $w$ in $B\big([0,1)\big)$ for  finitely many times.
 Since $h$ is a conformal map from $\mathbb{D}$
onto $\mathbb{D}-[0,1) $ and $\phi=B\circ h,$ then for each $w\in
\mathbb{D}-\{0\},$    $\phi-w$ is never a finite Blaschke product.
\end{proof}
\begin{rem}
In the proof of Lemma \ref{00l}, $\phi=B\circ h,$  where $B=\varphi_a B_0 .$ If we replace a finite Blaschke
product $B_1$ with $\varphi_a$ $(\,\varphi_a(z)=\frac{a-z}{1-\overline{a}z}\,)$, such that $B_1$ has no zeros on $[0,1)$, and that the restriction of
$\mathrm{arg} (B_1 )$ on $[0,1] $ attains any value in its range for
less than $n$ times for some positive integer $n$, then Lemma \ref{00l} also holds. For example,
 define
$$B_1 =\prod_{1\leq j\leq k}\varphi_{a_j},$$ where all $a_j$ lie in $\{w\in \mathbb{D}:
\frac{\pi}{2} <\arg w<\pi\}.$
Then $B_1$ has the desired property because for each $a_j$, $\arg \varphi_{a_j}$ can be
defined to be a continuous strictly-increasing function on $[0,1].$

 Rewrite $B=B_0B_1$, and put $\phi=B\circ h.$ As one will see later, the inner part of \mbox{$\phi-w$}  $(w\in
\mathbb{D})$ is a finite Blaschke product if and only if $w=0$.
Thus, the inner part of $\phi-\phi(\lambda)$ $(\lambda \in \mathbb{D})$ is a finite Blaschke product  if and only if $\lambda$ lies in
the zero set $Z(B_1)$ of $B_1,$ a finite subset of $\mathbb{D}$.

\end{rem}
 \noindent \textbf{Proof of Theorem \ref{351}:}
 The difficulty lies in the remaining part. That is to show  the inner part of $\phi$ is a finite Blashcke product.
 For this,  first we will give two computational results.

Put $S_1(z)=\exp(-\frac{1+z}{1-z})$. Clearly, it is continuous on
the unit circle except for $z=1$. We will see that $S_1(z)$
\emph{has non-tangential limit $0$ a}t $z=1.$ That is to show, for
each $\theta_0$ with $0<\theta_0<\frac{\pi}{2}$,
\begin{equation}\lim_{ \varepsilon \to 0^+,  |\theta|\leq \theta_0}S_1(1- \varepsilon e^{i \theta })=0. \label{24}\end{equation}
For this, write $z=1- \varepsilon e^{i \theta }$, where
$\varepsilon(\varepsilon>0)$ is enough small such that $z\in
\mathbb{D}$ whenever $|\theta|\leq \theta_0$. By direct
computations,
$$\mathrm{Re } \big(-\frac{1+z}{1-z}\big)=-\mathrm{Re } \big[ \frac{(1+z)(1-\overline{z})}{|1-z|^2}\big]
=-\frac{2\varepsilon \cos \theta-\varepsilon^2}{ \varepsilon^2},$$
which shows that
$$|S_1(1- \varepsilon e^{i \theta })|=\exp\big(-\frac{2\varepsilon \cos \theta-\varepsilon^2}{ \varepsilon^2}\big)
\leq \exp \big(-\frac{2  \cos \theta_0 }{ \varepsilon }+1 \big)\to
0, \, (\varepsilon\to 0^+).$$ Thus $S_1(z)$ has non-tangential limit
$0$ at $z=1.$ By the same computations, it follows that for any
$t>0$, $S_1^t(z) \triangleq \exp(-t\frac{1+z}{1-z})$ also has the
non-tangential limit $0$ at $z=1.$

Next another estimate will be given for $B$. Take a $\theta_1$ with
$0<|\theta_1|<\frac{\pi}{2}$. For example, take $\theta_1=\pm
\frac{\pi}{4}$. Then we have
\begin{equation}
\liminf_{m\to\infty} |B(1-\frac{1}{m!}e^{i\theta_1}) |>0.
\label{25}
\end{equation}
As before, let $d $ denote the pseudohyperbolic metric defined on
$\mathbb{D}$. Since for an enough large integer $m $,
$1-\frac{1}{m!} e^{i\theta_1}\in \mathbb{D}$. Then
\begin{eqnarray*}d(1-\frac{1}{n!},       1-\frac{1}{m!} e^{i\theta_1})& = &
\big|  \frac{  1-\frac{1}{n!}-(1- \frac{1}{m!}e^{i\theta_1})   }{   1-( 1-\frac{1}{n!})(1- \frac{1}{m!}e^{-i\theta_1}) } \big|  \\
 & = & \big| \frac{   \frac{1}{n!}- \frac{1}{m!}e^{i\theta_1}    }{   \frac{1}{n!}+ (\frac{1}{m!}-\frac{1}{n!m!})e^{-i\theta_1} } \big|,
 \end{eqnarray*}
and thus,
\begin{eqnarray} \prod_{n>m} d(1-\frac{1}{n!}, 1-\frac{1}{m!} e^{i\theta_1}) & \geq  &
 \prod_{n>m}  \dfrac{ |  \frac{1}{n!}- \frac{1}{m!}|    }{   \frac{1}{n!}+  \frac{1}{m!} } \nonumber \\
 & \geq  &    \dfrac{  1- \frac{1}{m+1}   }{ 1+  \frac{1}{m+1} }
  \prod_{n\geq m+2} \dfrac{   \frac{1}{m!}- \frac{1}{n!}   }{   \frac{1}{m!}+  \frac{1}{n!} } \nonumber  \\
  & \geq  &  \frac{1-\frac{1}{3}}{1+\frac{1}{3}}\prod_{k=2}^\infty \dfrac{ 1- \frac{1}{k!}  }{  1+  \frac{1}{k!} }  \nonumber  \\
   & > & \frac{1}{2} \prod_{k=2}^\infty \dfrac{ 1- \frac{1}{k!}  }{  1+  \frac{1}{k!} }\equiv \frac{1}{2}c>0. \label{26}
 \end{eqnarray}
Similarly, we have
\begin{eqnarray} \prod_{n<m} d(1-\frac{1}{n!}, 1-\frac{1}{m!} e^{i\theta_1}) & \geq  &
 \prod_{n<m}  \dfrac{   \frac{1}{n!}- \frac{1}{m!}  }{   \frac{1}{n!}+  \frac{1}{m!} } \nonumber \\
 & \geq  &  \frac{1-\frac{1}{3}}{1+\frac{1}{3}} \prod_{k=2}^{m-1} \dfrac{ 1- \frac{1}{k!}  }{  1+  \frac{1}{k!} }\nonumber \\
  & \geq  &  \frac{1}{2} \prod_{k=2}^\infty \dfrac{ 1- \frac{1}{k!}  }{  1+  \frac{1}{k!} } \nonumber \\
 & \geq  &   \frac{1}{2}  c>0. \label{3.10}
 \end{eqnarray}
Also, it is easy to see that $$\lim_{n\to\infty}d(1-\frac{1}{n!},
1-\frac{1}{n!} e^{i\theta_1})= \frac{  | 1-  e^{i\theta_1} |}{
|1+e^{i\theta_1} |}> 0. $$ Combining the above identity along with
(\ref{26}) and  (\ref{3.10}) shows that
\begin{equation} \liminf_{m\to\infty} |B(1-\frac{1}{m!}e^{i\theta_1}) |>0, \label{039}
\end{equation} as desired.
\vskip2mm
    \begin{figure} 
   \centering
     \includegraphics[scale=0.4]{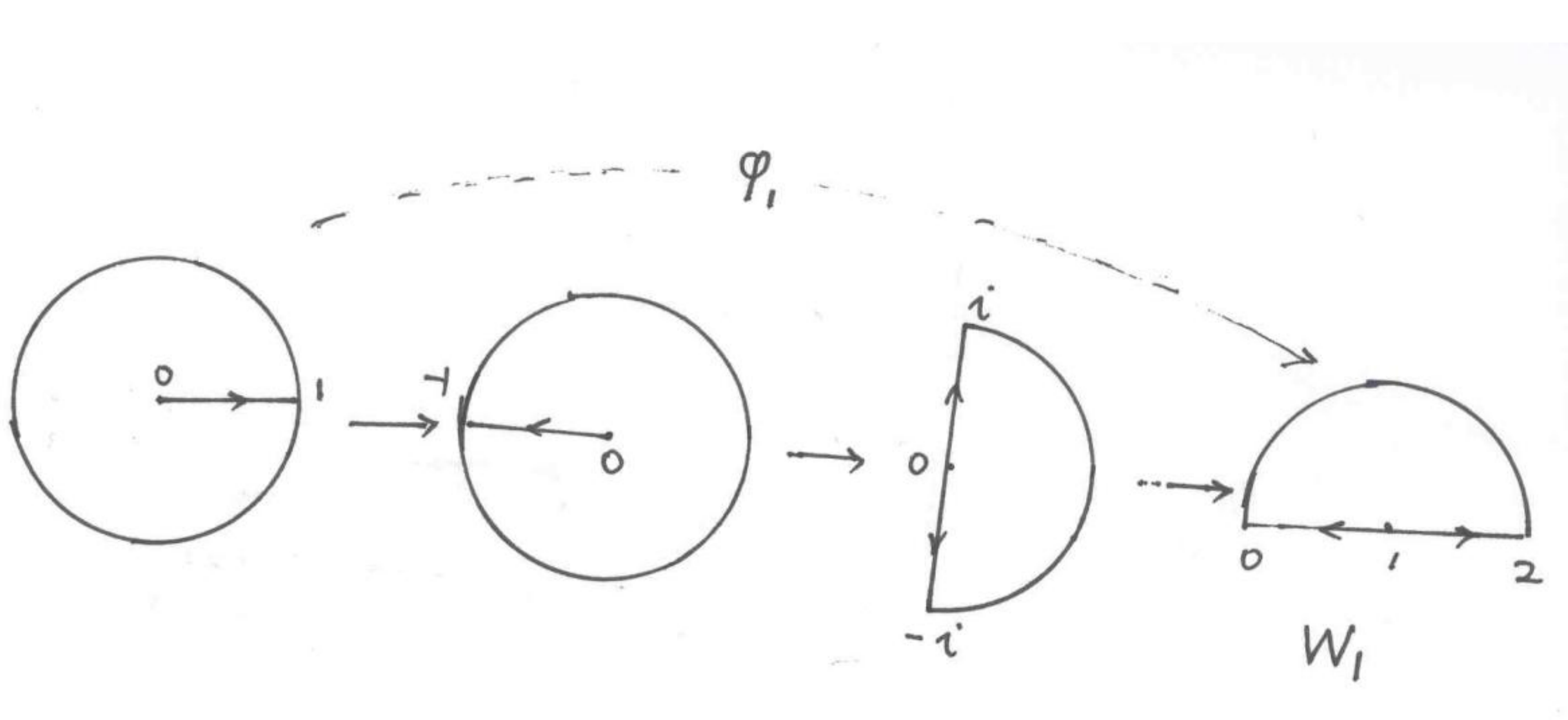} 
   \caption{$\varphi_1$} \label{feig3}
   \end{figure}
     \begin{figure}       
   \centering
     \includegraphics[scale=0.5]{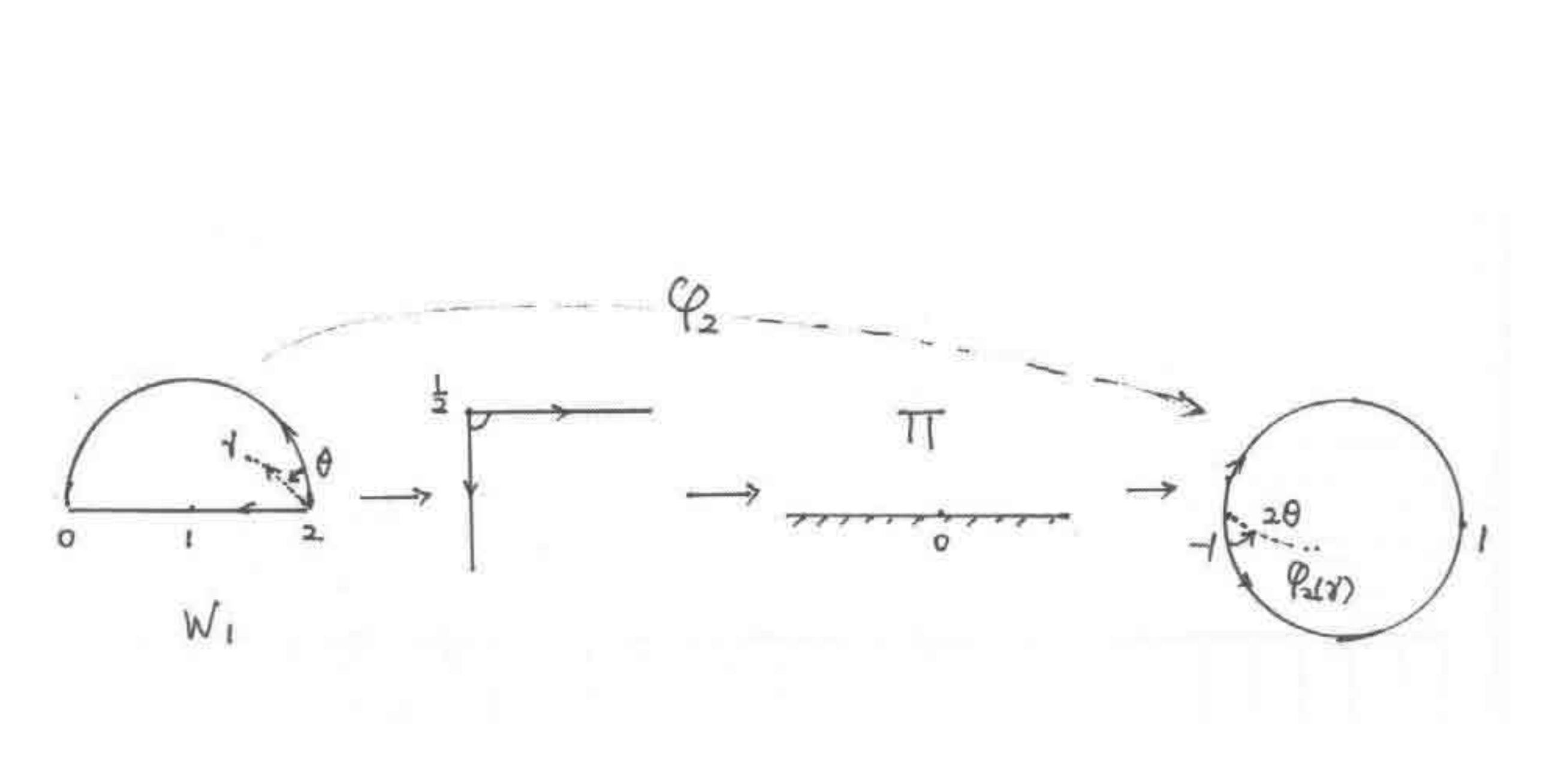} 
   \caption{$\varphi_2$} \label{feig2}
   \end{figure}
The idea is to compare (\ref{24}) with (\ref{039}) to derive a
contradiction. Below we shall give the concrete construction of
$h=g^{-1}$. Precisely, we will construct two conformal maps
$\varphi_2$ and $\varphi_1$, and $g \triangleq \varphi_2\circ
\varphi_1.$ Now define
$$\varphi_1(z)=i \sqrt{-z}+1, \, z\in \mathbb{D}-[0,1),$$ where
$\sqrt{1}=1,$ see Figure \ref{feig3}. Let us have a look at the
geometric property of $\varphi_1.$ Observe that $z\mapsto -z$ is a
rotation which maps $\mathbb{D}-[0,1) $ conformally to
$\mathbb{D}-(-1,0]  $.  A map is conformal if it preserves the angle
between two differentiable arcs. The  map $z\mapsto \sqrt{-z}$ is
conformal on $\mathbb{ C}- [0,+\infty)$, and hence on \linebreak
$\mathbb{D}-[0,1) $. One may think that the segment $[0,1]$ is split
into the upper and down parts, which are mapped onto two segments:
$i[-1,0] $ and $i[0,1]$. In particular, the point $1$ is mapped to
two points: $-i $ and $i$.
 Then with a rotation and a translation, $\varphi_1$ maps $\mathbb{D}-[0,1)$ conformally to the upper half disk
 $$W_1\triangleq \{\mathrm{Re}\, z>0; |z-1|<1\}.$$
Now  a biholomorphic map $\varphi_2:W_1\to \mathbb{D}$ will be
constructed. First, $z\mapsto \frac{1}{z}$ maps the the upper
half-disk $W_1$ onto a   rectangular domain $W_2$ between two half
lines: $\{\frac{1}{2}+it:t<0 \}$ and $[\frac{1}{2},+\infty)$. Then
with a rotation and a translation, $W_2$ is mapped onto the first
quadrant $W_3$. Write  $\varphi_3 (z)= z^2$, which maps  $W_3$ onto
the upper half plane $\prod$. Then one can give a mapping which maps
$\prod$ onto the unit disk, say, $z\mapsto \frac{z-i}{z+i}.$ Define
$\varphi_2$ to be the composition of the above maps, and we have
$$\varphi_2(z)=\frac{(\frac{1}{z}-\frac{1}{2})^2+i}{(\frac{1}{z}-\frac{1}{2})^2-i},\, z\in W_1.$$
Some words are in order. All the above maps are conformal; and
except for the map $\varphi_3: z\to z^2\ (z\in W_3),$ all maps are
conformal at the boundaries  of their domains of definition.
However, if $\theta(|\theta|<\frac{\pi}{2})$ is the angle between
two differential arcs beginning at $z=0$,
  then the angle between their image-arcs under $\varphi_3$ equals $2\theta$, see Figure \ref{feig2}.
Since $h=g^{-1}$  and $g = \varphi_2\circ
\varphi_1,$ by some computations  $$h(z)=\Big(
\frac{2\lambda(z)-1}{2\lambda(z)+1} \Big)^2$$where$$\lambda(z)=\sqrt{-i
\frac{z+1}{z+1}},$$ where $\sqrt{\cdot }$ denote the branch defined
on $\mathbb{C}-[0,+\infty)$ satisfying $\sqrt{-1}=-i.$

 After some verification, one sees that
the   map $h :\mathbb{D}\to \mathbb{D}-[0,1)$ extends
continuously onto $\overline{\mathbb{D}}$, which maps exactly two
points $ \eta_1$ and $\eta_2$ on $\mathbb{T}$ to $1$; $\mathbb{T}$
onto $\partial (\mathbb{D}-[0,1))$; one arc
$\widetilde{\eta_1\eta_2}$   onto $[0,1)$ for twice. In fact, by
some computations we have $\eta_1=1$ and $\eta_2=-1.$ Notice that
any non-tangential domain at $\eta_1$ or $\eta_2$ will be mapped to
some non-tangential domain at $1$, lying either above or below the
real axis, and vice versa. By the latter term ``non-tangential", we
mean the boundary of domain is not tangent to   $\mathbb{T}$  nor to
the segment  $[0,1]$ at $1$.

Now let $S$ denote the inner part of $\phi=B\circ h$. Observe that
$h^{-1}(0)$ contains exactly one point on $\mathbb{T}$, say
$\eta_0=h^{-1}(0)$. Since  $B$ is holomorphic on
$\overline{\mathbb{D}}-\{1\}$ and $h$ is holomorphic on
$\overline{\mathbb{D}}$ except for three possible points:
$$h^{-1}\{0,1\}=\{\eta_0, \eta_1, \eta_2\},$$
 one sees that $\phi=B\circ h$    is holomorphic at
 any point $\zeta\in \mathbb{T}-\{\eta_0, \eta_1,\eta_2\} $, and hence so is  $S$ \cite{Hof}.

 Next one will see that none of $\eta_0, \eta_1$ and $\eta_2$ is a singularity. Let
 $\phi=SF$ be the inner-outer decomposition of $\phi$, and then
 $|F|=|\phi|$, a.e. on $\mathbb{T}$. This shows that $F$ is bounded on $\mathbb{D}$\cite{Hof}.
It is easy to check that  $\phi $ is continuous at $\eta_0$,  and
$\phi(\eta_0)=B(0)\neq 0,$
   by  $\phi=SF$ one sees that $$\liminf_{z\to \eta_0}|S(z )|> 0.$$ Then  \cite[p.80, Theorem 6.6]{Gar},
   $\eta_0$ is not a singularity of $S$.
   Therefore, $\eta_1 $ and $\eta_2$ are the only possible singularities of $S$, and thus
  the singular part of $S$ is supported on $\{\eta_1,\eta_2\}$.
Remind that $S_1(z)=\exp(-\frac{1+z}{1-z})$. Put $b=h^{-1}(a),$ and
write
$$S(z)= \varphi_{b}(z) S_1^{t_1}(\overline{\eta_1}z)S_1^{t_2}(\overline{\eta_2}z),$$
where $t_1, t_2\geq 0 $. We will show that $t_1=t_2=0$ to finish the
proof.

For this, assume conversely that  either $t_1\neq 0$ or $t_2\neq 0$.
Without loss of generality, $t_1\neq 0.$  Then    $S $ has
non-tangential limit $0$ at $\eta_1$, and so does $\phi=SF$, where
$F$   is bounded on $\mathbb{D}$. However, with $\theta_1=\pm
\frac{\pi}{4}$,  and put
  $$\{z_k^1\}= \big{\{}h^{-1}(1-\frac{1}{k!} e^{\frac{\pi}{4}i} )\big{\}} \quad \mathrm{and} \quad \{z_k^2\}=
   \big{\{}h^{-1}(1-\frac{1}{k!} e^{- \frac{\pi}{4}i})\big{\}},$$
  where $k\geq n_0$ for some enough large integer $n_0 $ such that both $\{1-\frac{1}{k!} e^{\frac{\pi}{4}i}\} $ and
   $\{1-\frac{1}{k!} e^{- \frac{\pi}{4}i}\}$ lie in $\mathbb{D}$.
Since $\phi=B\circ h,$ by (\ref{25}) we get
\begin{equation} \liminf_{k\to\infty} |\phi(z_k^j)|> 0, \ \, j=1,2. \label{cont}\end{equation}
Considering  $$h^{-1}(1)=\{\eta_1,\eta_2\},$$ one notices that
$\{z_k^1\}$  and  $\{z_k^2\}$ are two  non-tangential sequences, one
tending to $\eta_1$ and the other to $\eta_2$. By (\ref{cont}),
this is a contradiction to that $\phi$ has non-tangential limit $0$
at $\eta_1$.  Therefore, $t_1=t_2=0$, and hence the inner part   $S$
of $\phi$ is  a M\"obius map. The proof of Theorem \ref{351} is complete. $\hfill
\square$
 \vskip2mm

Also, one can present more examples. It is worthwhile to mention that the sequence $\{1-\frac{1}{n!}\}_{n\geq 2}$
can be replaced with any thin Blaschke sequence $\{1-\varepsilon_n\}$ in $[0,1]$.
The reasoning is as follows. Without loss of generality, let $\{\varepsilon_n\}$ be decreasing to $0.$
Notice that if $0<|\theta_1|<\frac{\pi}{2},$ then
\begin{eqnarray*}d(1- \varepsilon_n,      1-  \varepsilon_m e^{i\theta_1})& = &
\big|  \frac{  1-\varepsilon_n-(1- \varepsilon_m e^{i\theta_1})   }{   1-(1- \varepsilon_n)(1- \varepsilon_m e^{-i\theta_1}) } \big|  \\
 & = & \big| \frac{  \varepsilon_n- \varepsilon_m e^{i\theta_1}    }{  \varepsilon_n+ (\varepsilon_m-\varepsilon_n\varepsilon_m)e^{ i\theta_1} } \big|\\
  & \geq  & \big| \frac{   \varepsilon_n- \varepsilon_m}{  \varepsilon_n+ (\varepsilon_m-\varepsilon_n\varepsilon_m)  } \big| \\
 & = &   d(1- \varepsilon_n,      1-  \varepsilon_m ) .
 \end{eqnarray*}
Since $\{1-\varepsilon_n\}$ is a thin Blaschke sequence; that is, $$\lim_{m\to\infty} \prod_{n;n\neq m}d(1- \varepsilon_n,      1-  \varepsilon_m )=1,$$
then it follows that
$$\lim_{m\to\infty} \prod_{n;n\neq m}d(1- \varepsilon_n,      1-  \varepsilon_me^{i\theta_1} )=1,$$
forcing $$\liminf_{m\to\infty} |B(  1-  \varepsilon_me^{i\theta_1})|>0.$$
This is a generalization of (\ref{25}), and the next discussion is just the same.
 
 \vskip3mm \noindent{Kunyu Guo, School of Mathematical Sciences,
Fudan University, Shanghai, 200433, China, E-mail:
kyguo@fudan.edu.cn}

\noindent{Hansong Huang, Department of Mathematics,
East China University of Science and Technology,  Shanghai, 200237, China, E-mail: hshuang@ecust.edu.cn
\end{document}